%% file: main.tex
\documentclass[letterpaper, 10 pt, conference]{ieeeconf}
\IEEEoverridecommandlockouts                              
\usepackage{soul}
\usepackage{comment}
\usepackage{graphics}
\usepackage{epsfig}
\usepackage{amsmath}
\usepackage{amssymb}
\usepackage{cite}

\usepackage{bm,upgreek}
\newtheorem{theorem}{Theorem}

\newtheorem{definition}{Definition}
\newtheorem{remark}{Remark}
\newtheorem{assumption}{Assumption}
\usepackage{xcolor}
\usepackage{float}

\input{formulas}
\usepackage{color}
\definecolor{mypurple}{rgb}{0.4, 0.06, 0.9}
\definecolor{mymaroon}{rgb}{0.8, 0.1, 0.5}
\definecolor{mydarkblue}{rgb}{0.0, 0.06, 0.8}
\usepackage{hyperref}
\hypersetup{
    colorlinks=true,  
    linkcolor=mymaroon,  
    citecolor=mydarkblue,   
    urlcolor=mypurple,    
    filecolor=blue,   
    pdfborder={0 0 0} 
}

\title{\LARGE \bf Safe Control for Pursuit-Evasion with Density Functions}

\author{Mustafa Bozdag, Arya Honarpisheh, Mario Sznaier
\thanks{\textcolor{black}{This work was partially
supported by
NSF grants CNS--2038493 and CMMI--2208182, AFOSR grant FA9550-19-1-0005, ONR grant N00014-21-1-2431, and  the Sentry DHS Center of Excellence
under Award 22STESE00001-03-03.}
 The authors are with the Robust Systems Lab, ECE Department, Northeastern University, Boston, MA 02115. (e-mails: {bozdag.m,honarpisheh.a}@northeastern.edu,
msznaier@coe.neu.edu)
}}%

\begin{document}

\maketitle
\thispagestyle{empty}
\pagestyle{empty}

\begin{abstract}

This letter presents a density function based safe control synthesis framework for the pursuit-evasion problem. We extend safety analysis to dynamic unsafe sets by formulating a reach-avoid type pursuit-evasion differential game as a robust safe control problem. Using density functions and semi-algebraic set definitions, we derive sufficient conditions for weak eventuality and evasion, reformulating the problem into a convex sum-of-squares program solvable via standard semidefinite programming solvers. This approach avoids the computational complexity of solving the Hamilton-Jacobi-Isaacs partial differential equation, offering a scalable and efficient framework. Numerical simulations demonstrate the efficacy of the proposed method.

\end{abstract}

\section{INTRODUCTION}

Safety plays a critical role in control theory applications as the systems we deploy in the real world must operate reliably. This requirement has become as fundamental as classical notions like stability or controllability, especially with the widespread use of cyber-physical systems such as drone swarms and autonomous vehicles \cite{alan2023autonomous, buyukkocak2024sequential}. Designing controllers for such systems is often challenging as they must satisfy multiple objectives and constraints to ensure safety.

Numerous methods have been developed to certify the safety of dynamical systems since Nagumo introduced set invariance in the 1940s \cite{nagumo1942invarient}. Set invariance has been central to many of these approaches following the advent of Lyapunov theory. Level-set methods separate regions of interest using the 0-contour of a function, with barrier functions ensuring trajectory safety if their superlevel sets remain invariant \cite{prajna2004hybrid}. Control barrier functions (CBFs) extend this concept, akin to the generalization of Lyapunov functions to control Lyapunov functions (CLFs) \cite{sontag1983lyapunov}, enabling the synthesis of safe control laws \cite{ames2017quadratic}. Density functions, originally formulated as a dual to Lyapunov stability, serve as barrier-like certificates \cite{rantzer2001dual_density_lyapunov, rantzer2004analysis_synthesis_safety}.

Like barrier functions, they define safe and unsafe regions via a 0-level set, ensuring safety in the superlevel sets under invariance conditions. Safe controllers and density functions can be jointly computed using convex methods, including sum-of-squares (SOS) formulations \cite{prajna2004nonlinear, yu2021convex}. Data-driven approaches \cite{yu2022data} and robust extensions \cite{zheng2023density} further enhance these methods. The field of safe control primarily addresses problems with static unsafe sets. Static reachable sets extend safety analysis to reach-avoid scenarios, often characterized by Lyapunov-like reachability conditions \cite{xue2024reach, xue2024synthesis, zheng2024ddsc}. The reach-avoid problem is also widely studied in the Hamilton-Jacobi (HJ) reachability framework from optimal control theory \cite{bansal2017hamilton, bryson2018applied_optimal_control}. The HJ equation, a partial differential equation (PDE) governing the value function’s evolution, plays a key role in differential game theory alongside Pontryagin's maximum principle.

Pursuit-evasion games are zero-sum differential games extensively studied in \cite{isaacs1999differential}. A classical two-player pursuit-evasion game is generally formulated as a minimax problem with time-to-capture as the objective. In the reach-avoid case, reaching a target for the evader flags the termination of the game instead of a time constraint \cite{weintraub2020pursuit}. There is a growing body of work in the literature addressing alternative settings, including perimeter-defense games, where geometric solution methods are employed instead of directly solving HJI equations \cite{shishika2020cooperative, bajaj2024multivehicle, shishika2021partial}. Although geometric methods could be applied to solve such a problem, under more complex state constraints and dynamics, these methods become impractical \cite{oyler2016pursuit, bera2017comprehensive_game_of_two_cars}. In \cite{margellos2011hamilton}, the value function for a pursuit-evasion game is formulated and shown to be the unique viscosity solution to the Hamilton-Jacobi-Isaacs (HJI) equation. 
However, solving the HJI PDE subject to state and time constraints is a daunting task. To solve an HJI equation on the reachable sets, discretization of the state space is often required, resulting in an exponential scaling of computational complexity with respect to system dimensionality.

In this work, we formulate the pursuit-evasion game as a robust safe control problem with a dynamic unsafe set, where the evader must reach a target set while avoiding capture by a pursuer and staying within a bounded environment. The unsafe region is defined by a catch radius around the pursuer. By modeling both agents as a single system and treating the pursuer's action as a bounded disturbance, we recast the problem as robust control synthesis. Using density functions, we derive sufficient conditions for safety and reachability. Assuming polynomial dynamics and semi-algebraic sets, we apply Putinar’s Positivstellensatz to convert the problem into a sum-of-squares (SOS) program solvable via semidefinite (SDP) programming. To the best of our knowledge, our work is the first study of pursuit-evasion games through the lens of safe control methods in a computationally efficient fashion. The main contributions of this work consist of:
\begin{itemize}
    \item Formulating the pursuit-evasion game as a robust safe control problem with the use of density functions,
    \item Proposing a convex optimization framework as an alternative to the numerically demanding task of solving PDEs with state constraints,
    \item Providing a numerical simulation to validate the theoretical approach, demonstrating its practical feasibility.
\end{itemize}

\section{PRELIMINARIES}
\label{sec:prelim}
\subsection{Notation}\label{Notation}
\vspace{-18pt}
\begin{alignat*}{2}
    & \mathbb{R},\mathbb{R}^n \quad&&  \text{Set of  real numbers, \textit{n}-tuples } \\
    & x,\mathbf{x},\mathbf{X}  \quad&& \text{Scalar, vector, matrix} \\
    & \textbf{1}, \textbf{0}, \mathbf{I} \quad && \text{Vector/matrix of all 1s, 0s, identity matrix} \\
    & \|\mathbf{x}\|_{\infty} && \ell_{\infty} \text{-norm of vector $\mathbf{x}$} \\
    & \|\mathbf{x}\|_{2} && \ell_{2} \text{-norm of vector $\mathbf{x}$} \\
    & \mathbf{X} \succeq 0 && \mathbf{X} \text{ is positive semidefinite} \\
    & f \in \mathcal{C}^d && \text{The $d^{th}$ partial derivatives of $f$ are continuous} \\
    & \nabla V && \text{The gradient of a scalar function $V$} \\
    & \nabla \cdot f && \text{The divergence of a vector valued function $f$}\\
    & \mathbb{R}[\mathbf{x}]& & \text{Polynomials in the indeterminate $\mathbf{x} \in \mathbb{R}^n$}
\end{alignat*}

\subsection{Safety and Eventuality}
The general formulation of a nonlinear system is given as:
\begin{equation}
\dot{\mathbf{x}}=f(\mathbf{x},\mathbf{u}(\mathbf{x}))
    \label{eq:vanilla_dynamics}
\end{equation}
where $\mathbf{x}\in\mathcal{X}\subseteq\mathbb{R}^n$ and $\mathbf{u}: \mathcal{X} \rightarrow \mathcal{U}\subseteq \mathbb{R}^m$ are the state and the control input of the system respectively. Given the initial set $\mathcal{X}_i$, the unsafe set $\mathcal{X}_a$, and target set $\mathcal{X}_r$, we assume they are mutually exclusive and contained in the bounded set $\mathcal{X} \subset \mathbb{R}^n$. Assuming $f$ and $\mathbf{u}$ belong to $\mathcal{C}^1$, and given the initial condition $\mathbf{x}(0) = \mathbf{x}_0 \in \mathcal{X}_i$, the solution to \eqref{eq:vanilla_dynamics} exists and is unique. Let $\phi_t(\mathbf{x}_0)$ denote this solution as a function of time $t$. The system \eqref{eq:vanilla_dynamics}:

\begin{enumerate}
\item Is \emph{safe under $\mathcal{X} \setminus \mathcal{X}_a$} if there exists $\mathbf{u}: \mathcal{X} \rightarrow \mathcal{U}$ such that for all $ \mathbf{x}_0 \in \mathcal{X}_i$, $\forall t \in [0,\infty), \; \phi_t(\mathbf{x}_0) \in \mathcal{X} \setminus \mathcal{X}_a$;

\item \emph{Eventually reaches $\mathcal{X}_r$} if there exists $\mathbf{u}: \mathcal{X} \rightarrow \mathcal{U}$ such that for all $\mathbf{x}_0$ in $\mathcal{X}_i$, $\exists T \in [0,\infty), \; \phi_T(\mathbf{x}_0) \in \mathcal{X}_r$;
 
\item Is \emph{safe under $\mathcal{X} \setminus \mathcal{X}_a$ until it eventually reaches $\mathcal{X}_r$} if there exists $\mathbf{u}: \mathcal{X} \rightarrow \mathcal{U}$ such that for all $\mathbf{x}_0$ in $\mathcal{X}_i$, $\exists T \in [0,\infty), \; \phi_T(\mathbf{x}_0) \in \mathcal{X}_r$ and $\forall t \in [0,T] \; \phi_t(\mathbf{x}_0) \in \mathcal{X} \setminus \mathcal{X}_a$.
\end{enumerate}
If the eventuality property holds for almost everywhere in $\mathcal{X}_i$, then the \emph{weak eventuality} holds \cite{aubin2011viability}. The term ``almost everywhere” (and ``almost all” later on) implies the property holds at all states except on a set of measure zero.

\subsection{Pursuit-Evasion Games}

The core problem in differential game theory is the zero-sum two-player game known as pursuit-evasion for the following dynamical system:
\begin{equation}
    \begin{aligned}
    \dot{\mathbf{x}}=f(\mathbf{x},\mathbf{u}(\mathbf{x}),\mathbf{w}(\mathbf{x}))
    \end{aligned}
    \label{eq:vanilla_pe_dynamics}
\end{equation}
where the additional $\mathbf{w}: \mathcal{X} \rightarrow \mathcal{W}\subseteq \mathbb{R}^m$ is the control action of the pursuer. In the vanilla case, the goal of the pursuer is to catch the evader while the evaders goal is to avoid being captured as long as possible, with time-to-capture as the cost/reward for the agents \cite{isaacs1999differential, weintraub2020pursuit}. We focus on the reach-avoid case of the game from the evader's perspective, where the evader is trying to \textbf{reach} a target set $\mathcal{X}_r$, rather than  maximizing time-to-capture, while \textbf{avoiding} an unsafe set $\mathcal{X}_a$, formulated as the following:
\begin{equation}
    \begin{aligned}
        \mathcal{X}_r = \{ \mathbf{x}: (h_r(\mathbf{x}) \leq 0)\}, \; \mathcal{X}_a = \{ \mathbf{x}: (h_a(\mathbf{x}) \leq 0)\}.
    \end{aligned}
    \label{eq:reach_avoid_sets}
\end{equation}
where $h_r$ and $h_a$ are functions whose zero-level sets implicitly define the boundaries of the reach and avoid regions. The traditional approach to solve this problem involves the construction of a value function:
\begin{equation}
    \begin{aligned}
        V(\mathbf{x},t)=\inf_{\mathbf{u}}{\sup_{\mathbf{w}}{\max{\left\{ h_r(\phi_T), \; -\min_{\tau\in[t,T]}{h_a(\phi_{\tau})} \right\}}}}
    \end{aligned}
    \label{eq:value_func}
\end{equation}
which is the unique viscosity solution of the following Hamilton-Jacobi-Isaacs (HJI) PDE \cite{margellos2011hamilton}:
\begin{equation}
    \begin{aligned}
        \max \left\{ -h_a(\mathbf{x})-V(\mathbf{x},t), \; \frac{\partial V}{\partial t} + \sup_{\mathbf{u}}{\inf_{\mathbf{w}}{\frac{\partial V}{\partial \mathbf{x}}\dot{\mathbf{x}}}} \right\} = 0
    \end{aligned}
    \label{eq:hji}
\end{equation}
Solving the HJI equation for reach-avoid pursuit-evasion games is inherently challenging due to the safety and reachability constraints, and the need to handle min-max optimization over the control strategies of both players. Instead, we try to approach this problem from a safe control perspective.

\subsection{Density Functions}

Density functions, introduced in \cite{rantzer2001dual_density_lyapunov} as a dual to Lyapunov stability, represent a distribution over states evolving according to the system dynamics. The following theorem provides sufficient conditions for the weak eventuality and safety properties using density functions for autonomous systems.

\begin{theorem}
\label{thm:weak-eventuality-rantzer}
From \cite{prajna2007convex}. Consider the dynamical system $\mathbf{\dot{x}}=f(\mathbf{x})$ with $f\in \mathcal{C}^1(\mathbb{R}^n,\mathbb{R}^n)$, the bounded set $\mathcal{X}\subset\mathbb{R}^n$, and the sets $\mathcal{X}_i,\mathcal{X}_a,\mathcal{X}_r\subseteq \mathcal{X}$. If there exists an open set $\Tilde{\mathcal{X}_i}\supseteq\mathcal{X}_i$ and a density function $\rho\in \mathcal{C}^1(\mathbb{R}^n)$ satisfying:
   \begin{align}
    \rho(\mathbf{x}) &\geq 0, \; 
    \forall \mathbf{x} \in \Tilde{\mathcal{X}_i}, \label{eq:rantzer_initial_condition} \\
    \rho(\mathbf{x}) &< 0, \; 
    \forall \mathbf{x} \in \operatorname{cl}(\partial\mathcal{X} \setminus \partial\mathcal{X}_r) \cup \mathcal{X}_a,  \label{eq:rantzer_unsafe_condition} \\
    \nabla \cdot (\rho f)(\mathbf{x}) &> 0, \; \forall \mathbf{x} \in \operatorname{cl}(\mathcal{X} \setminus \mathcal{X}_r), \label{eq:rantzer_divergence_condition}
    \end{align}

then the \emph{weak eventuality and safety} properties hold. For almost all initial states $\mathbf{x}_0\in \mathcal{X}_i$, the trajectory $\mathbf{x}(t)$ with $\mathbf{x}(0)=\mathbf{x}_0$ will satisfy $\mathbf{x}(T)\in\mathcal{X}_r$ for some $T\geq 0$ and $\mathbf{x}(t)\notin \mathcal{X}_a$, $\mathbf{x}(t)\in \mathcal{X}$ for all $t\in [0,T]$.
\end{theorem}

The proof of Theorem \ref{thm:weak-eventuality-rantzer} relies on Liouville's Theorem given in \cite{rantzer2001dual_density_lyapunov}. Building on this result, in this letter, we extend the framework to pursuit-evasion games.
\begin{remark}The advantages of using a density based approach are two-fold: (i) the eventuality property guarantees that the evader will reach a safe target set and, (ii) 
as shown in Theorem \ref{thm:eventuality-safty}, 
it leads to a jointly convex formulation in the density function and evader strategy. Further, this strategy does not require computing the optimal pursuer policy.
\end{remark}

\subsection{Sum-of-Squares}

Sum-of-squares (SOS) optimization is a framework used in certain classes of polynomial optimization problems where the main goal is to establish non-negativity over a region. The SOS formulations in this paper are based on Putinar's Positivstellensatz, defined as the following \cite{parrilo2000sos,putinar1993psatz}.

\begin{definition}
A polynomial $ p \in \mathbb{R}[x] $ is SOS if there exist polynomials $ \{q_i \in \mathbb{R}[x]\}_{i=1}^N $ such that:
\begin{equation} \label{eq:sos_definition}
    p(x) = \sum_{i=1}^N q_i(x)^2.
\end{equation}
\end{definition}
The set of SOS polynomials forms the cone $ \Sigma[x] $, with its restriction to degree $ 2d $ denoted as $ \Sigma_d[x] $. The cone $ \Sigma_d[x] $ has a positive semidefinite (PSD) representation given by:
\begin{equation} \label{eq:sdp_representation}
    p(x) = v(x)^\top Q v(x),
\end{equation}
where $ v(x) $ is the vector of monomials up to degree $ d $, and $ Q \succeq 0 $ is a PSD matrix. A sufficient condition for $ p(x) $ to be nonnegative over the semialgebraic set:
\begin{equation} \label{eq:semialgebraic_set}
    \{x \mid h_i(x) \geq 0, \, i = 1, \dots, N_c\}
\end{equation}
is the existence of $ \sigma_0, \dots, \sigma_{N_c} \in \Sigma[x] $ such that: 
\begin{equation} \label{eq:sos_condition}
    p(x) = \sigma_0 + \sum_{i=1}^{N_c} \sigma_i h_i(x).
\end{equation}
In this work, we adopt a similar approach to \cite{zheng2023density, zheng2024ddsc}, leveraging SOS optimization for the simultaneous synthesis of density functions and control laws.

\section{Problem Statement}
\label{sec:problem_statement}
This paper aims to design a controller to guide an evader in a reach-avoid pursuit-evasion game. The dynamics of the pursuer and evader are described as follows:
\begin{equation}
\begin{aligned}
    \dot{\mathbf{x}}_e &= f_e(\mathbf{x}_e) + g_e(\mathbf{x}_e) \mathbf{u}(\mathbf{x}_e, \mathbf{x}_p) \\
    \dot{\mathbf{x}}_p &= f_p(\mathbf{x}_p) + g_p(\mathbf{x}_p) \mathbf{w}(\mathbf{x}_e, \mathbf{x}_p)
\end{aligned}
\label{eq:game_dynamics}
\end{equation}
where $\mathbf{x}_e \in \mathbb{R}^{n_e}$ and $\mathbf{x}_p \in \mathbb{R}^{n_p}$ represent the coordinates of the evader and pursuer, respectively. The variable $\mathbf{u} \in \mathbb{R}^{m_e}$ denotes the control input to the evader, while $\mathbf{w} \in \mathbb{R}^{m_p}$ denotes the input to the pursuer, which is considered uncontrollable by the evader. To represent both the evader and pursuer as a unified dynamical system, we define the combined state vector as:
\begin{equation}
\label{eq:combined_states_inputs}
\mathbf{x} = 
\begin{bmatrix} 
\mathbf{x}_e^\top & \mathbf{x}_p^\top 
\end{bmatrix}^\top.
\end{equation}
Therefore, the combined system is
\begin{equation}
\label{eq:combined_system}
    \mathbf{\dot{x}} = f(\mathbf{x}) + g_u(\mathbf{x}) \mathbf{u}(\mathbf{x}) + g_w(\mathbf{x}) \mathbf{w}(\mathbf{x})
\end{equation}
where
\begin{equation}
\label{eq:combined_dynamics}
    f = \begin{bmatrix}
        f_e \\ f_p
    \end{bmatrix},
    \quad
    g_u = \begin{bmatrix}
        g_e \\ \mathbf{0}
    \end{bmatrix},
    \quad
    g_w = \begin{bmatrix}
        \zero \\ g_p
    \end{bmatrix}.
\end{equation}
We assume that, initially, both the evader and pursuer are confined within the initial set $\mathcal{X}_i$. The objective of the evader is to reach a designated target set $\mathcal{X}_r$ while remaining within the bounded region $\mathcal{X}$ and avoiding capture by the pursuer. Capture occurs when the evader and pursuer are in proximity to each other, occupying a location within the unsafe set $\mathcal{X}_a$, which is a function of both $\mathbf{x}_e$ and $\mathbf{x}_p$. Additionally, we define the path-connected set $\mathcal{X}_c = \mathcal{X} \setminus \mathcal{X}_a$ as the region within which the evader can move freely.

\section{Density-Based Robust Evasion}
\label{sec:density_robust_control}
In pursuit-evasion, the evader typically does not have direct access to the pursuer policy. In principle, the optimal pursuer policy can be inferred by solving the PDE \eqref{eq:hji}, but, as mentioned before, this is far from trivial. In this paper, we will circumvent this difficulty by modeling the pursuer's strategy as an unknown but bounded disturbance $\mathbf{w}(\mathbf{x}) \in \mathcal{W}$ to be rejected by a robust control policy. In other words, regardless of the actions taken by the pursuer, the evader safely reaches the designated target. To this end, we need the following technical requirements: 
\begin{assumption}
\label{ass:assumption1}
Let the following conditions hold:
\begin{itemize}
    \item $\mathcal{X}$ is bounded,
    \item The initial set $\mathcal{X}_i$, unsafe set $\mathcal{X}_a$, and target set $\mathcal{X}_r$ are mutually disjoint subsets of $\mathbb{R}^n$ contained in $\mathcal{X}$,
    \item $\mathcal{X}_i$ is an open set,
    \item The active set $\mathcal{X}_c = \operatorname{cl} ( \mathcal{X} \setminus (\mathcal{X}_r \cup \mathcal{X}_a) )$ has a path-connected interior, i.e., at each time, there is a continuous path from the evader to the target within $\mathcal{X}_c$,
    \item $\mathcal{W}$ is a compact subset of $\mathbb{R}^{n_p}$.
\end{itemize}
\end{assumption}

\begin{theorem}
\label{thm:eventuality-safty}
    Consider the dynamical system \eqref{eq:combined_system} with $f \in \mathcal{C}^1(\mathbb{R}^n, \mathbb{R}^n)$, $g_u, g_w \in \mathcal{C}^1(\mathbb{R}^n, \mathbb{R}^m)$, $n=n_e+n_p$, and $m=m_e+m_p$. Let $V \in \mathcal{C}^1(\mathbb{R}^n, \mathbb{R}(.)$, $V(x)>0 \; \forall x \in \operatorname{cl} ( \mathcal{X} \setminus \mathcal{X}_r)$  and $\alpha \geq 1$ be a given function and scalar, respectively. Then,  if there exist a density function $\rhon: \mathbb{R}^n \rightarrow \mathbb{R}$ and a vector field $\psin: \mathbb{R}^n \rightarrow \mathbb{R}^n$ continuously differentiable on $\mathcal{X}_c$  such that Assumption \ref{ass:assumption1} and the following conditions hold
    \begin{align}
    &V(\mathbf{x}) > 0, \;
    \forall \mathbf{x} \in \mathcal{X} \setminus \mathcal{X}_r,\label{eq:V_condition} \\
    &\rhon(\mathbf{x}) \geq 0, \; 
    \forall \mathbf{x} \in \mathcal{X}_i, \label{eq:initial_condition} \\
    &\rhon(\mathbf{x}) < 0, \; 
    \forall \mathbf{x} \in \operatorname{cl}(\partial\mathcal{X} \setminus \partial\mathcal{X}_r) \cup \mathcal{X}_a,  \label{eq:unsafe_condition} \\
    &\nabla \cdot (\rhon (f + g_w \mathbf{w}) + g_u \psin) > \nonumber \\
    & 
    \quad \alpha \frac{\nabla V}{V} \cdot (\rhon (f + g_w \mathbf{w}) + g_u \psin),
    \label{eq:divergence_condition} \\
    & \qquad\qquad\;\; \forall \mathbf{w} \in \mathcal{W}, \; \forall \mathbf{x} \in \operatorname{cl}(\mathcal{X} \setminus \mathcal{X}_r), \nonumber
    \end{align}
    then, the control input $\mathbf{u} = \frac{\psin}{\rhon}$ enforces the \emph{weak eventuality and evasion} properties on \eqref{eq:combined_system} in a robust manner i.e. for almost all $\mathbf{x}_0 \in \mathcal{X}_i$, $\exists T \geq 0 \; \phi_T(\mathbf{x}_0) \in \mathcal{X}_r$ and $ \forall t \in [0,T] \; \phi_t(\mathbf{x}_0) \notin \mathcal{X}_a$ for all pursuers actions $\mathbf{w} \in \mathcal{W}$.
\end{theorem}

\begin{proof}
For the controlled system $\eqref{eq:combined_system}$, to apply the condition \eqref{eq:rantzer_divergence_condition} treating $\mathbf{w}$ as a disturbance, we must have $\nabla \cdot (\rho (f + g_u \mathbf{u} + g_w \mathbf{w})) > 0$ for all $\mathbf{w} \in \mathcal{W}$. To make this condition linear in the unknowns $\rho$ and $\bm{\uppsi}$, we introduce $\bm{\uppsi} = \rho \mathbf{u}$. The condition then becomes
$\nabla \cdot (\rho (f +  g_w \mathbf{w}) + g_u \bm{\uppsi}) > 0$ for all $\mathbf{w} \in \mathcal{W}$. Substitute $\rho = \frac{\rhon}{V^\alpha}$ and $\bm{\uppsi} = \frac{\psin}{V^\alpha}$ to express it in terms of $\rhon$ and $\psin$. Rewriting the LHS gives:
\begin{align*}
    & \nabla \cdot \left(\frac{\rhon}{V^\alpha} (f +  g_w \mathbf{w}) + g_u \frac{\psin}{V^\alpha} \right) \\
    &= \frac{V^\alpha \nabla \cdot (\rhon (f +  g_w \mathbf{w}) + g_u \psin)}{V^{2\alpha}} \\
    & \quad - \frac{\alpha V^{\alpha-1} \nabla V \cdot (\rhon (f +  g_w \mathbf{w}) + g_u \psin)}{V^{2\alpha}}.
\end{align*}
Since $V^\alpha > 0$, the condition $\nabla \cdot (\rho (f +  g_w \mathbf{w}) + g_u \bm{\uppsi} ) > 0$ holds if and only if:
\begin{equation*}
\nabla \cdot (\rhon (f +  g_w \mathbf{w}) + g_u \psin) - \alpha \frac{\nabla V}{V} \cdot (\rhon (f +  g_w \mathbf{w}) + g_u \psin) > 0,
\end{equation*}
which corresponds to \eqref{eq:divergence_condition}.
Because $V$ is positive on $\mathcal{X} \setminus \mathcal{X}_r$, the conditions \eqref{eq:initial_condition} and \eqref{eq:unsafe_condition} follow directly from \eqref{eq:rantzer_initial_condition} and \eqref{eq:rantzer_unsafe_condition}, respectively.
Thus, the three conditions \eqref{eq:initial_condition}, \eqref{eq:unsafe_condition}, and \eqref{eq:divergence_condition} imply the conditions of Theorem \ref{thm:weak-eventuality-rantzer}, completing the proof.
\end{proof}
Since Theorem \ref{thm:eventuality-safty} is framed for a pursuit-evasion game, we refer to the safety property as \emph{evasion}, highlighting its ability to handle dynamic unsafe sets. In Theorem~\ref{thm:eventuality-safty}, the feasibility conditions represent the \emph{winning regions} for the evader, and the controller $\mathbf{u} = \frac{\psin}{\rhon}$ is the \emph{security strategy} of the evader.
\begin{remark} Note that setting  $V(.) \equiv  1$ recovers the conditions in   Theorem \ref{thm:weak-eventuality-rantzer}. These additional degrees of freedom  facilitate the satisfaction of the eventuality condition via a rational density function, as in \cite{moyalan2024synthesizing, zheng2023control_density_lyapunov}. In particular, in this paper, motivated by \cite {zheng2023control_density_lyapunov} we use $V := \| \mathbf{x}_e - \mathbf{x}_r \|_2^2$, the squared distance to the center of the target arc $\mathbf{x}_r$.
\end{remark}
\begin{remark}
The controlled system \eqref{eq:combined_system} with $\mathbf{u} = \frac{\psin}{\rhon}$ becomes singular at points where $\rhon = 0$. To avoid these regions, it is necessary to introduce a constraint that bounds $\mathbf{u}$ along the safe trajectory, ensuring the system remains well-defined. Since $\rhon$ and $\psin$ are optimization variables, we must impose constraints to satisfy $\|\mathbf{u}\|_{\infty} = \|\frac{\psin}{\rhon}\|_{\infty} \leq u_{\max}$:  
\begin{equation}
    (\rhon u_{\max} - \psin \geq 0) \wedge (\rhon u_{\max} + \psin \geq 0).
\end{equation}  
 \end{remark}
\vspace{5pt}

\begin{figure}[ht]
    \centering
    \includegraphics[width=0.7\linewidth]{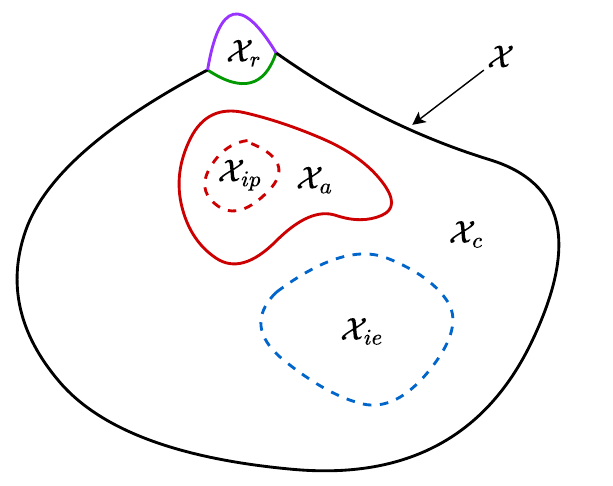}
    \caption{The reach-avoid environment. The objective is to remain in $\mathcal{X}$, avoiding $\mathcal{X}_a$, and reaching $\mathcal{X}_r$ using a connected path within $\mathcal{X}_c$. The initial states for the evader and the pursuer are indicated with $\mathcal{X}_{ie}$ and $\mathcal{X}_{ip}$, which together form $\mathcal{X}_{i}=\mathcal{X}_{ie}\cap\mathcal{X}_{ip}$. For pursuit-evasion, the unsafe set $\mathcal{X}_a$ is dynamic. It is important to note that all these sets are defined in a \textbf{four-dimensional state space}. The 2D illustration provided is for intuition only and does not fully capture the true structure of these sets.}
    \label{fig:environment}
    \vspace{-10pt}
\end{figure}

Figure \ref{fig:environment} illustrates the scenario of Theorem \ref{thm:eventuality-safty} from the evader’s perspective, simplifying the four-dimensional state-space to 2D. Notably, the purple curve does not influence the conditions of Theorem \ref{thm:eventuality-safty}. This curve represents the outer boundary of the the target set $\mathcal{X}_r$ which only appears in conditions \eqref{eq:unsafe_condition} and \eqref{eq:divergence_condition}. In these conditions, the relevant boundary of $\mathcal{X}_r$ is the inner green curve, because it forms a shared boundary between $\mathcal{X}_r$ and $\mathcal{X}_c$, implying that any trajectory within $\mathcal{X}_c$ that reaches $\mathcal{X}_r$ must necessarily cross it. Thus, $\mathcal{X}_r$ is only accessible from the green curve for trajectories within $\mathcal{X}_c$, rendering the purple curve irrelevant.

In standard reach-avoid scenarios,
the system under consideration cannot admit a stable equilibrium point in $\mathcal{X}_c$; otherwise, the divergence condition \eqref{eq:divergence_condition} would be violated. Thus, stable equilibria are typically excluded from $\mathcal{X}_c$ explicitly \cite{rantzer2004analysis_synthesis_safety}. In a pursuit-evasion game, this scenario is inherently precluded. Suppose, by contradiction, that the strategies $\mathbf{u(x)},\mathbf{w(x)}$ are such that the closed loop  system \eqref{eq:combined_system} has an equilibrium point $(\mathbf{x}_e^*, \mathbf{x}_p^*)$. For this equilibrium to exist, the pursuer must have no incentive to move. However, since there exists an unsafe region $\mathcal{X}_a$ that both agents cannot occupy simultaneously, $\mathbf{x}_e^*$ is necessarily outside of this region. If the evader is not within the unsafe region, the pursuer must move towards it to fulfill its objective, thereby contradicting the assumption that $\mathbf{x}_p^*$ is an equilibrium.

\section{Convex Relaxation}
\label{sec:convex_relaxation}
The density-based robust control formulation presented in Section~\ref{sec:density_robust_control} results in a nonlinear feasibility problem. To solve this problem efficiently, we make the following simplifying assumptions: (i) the dynamics are integrators with a constraint that bounds the velocity; (ii) the initial, target, and unsafe sets are semi-algebraic; and (iii) $\rho$ and $\bm{\uppsi}$ are polynomials. Under these assumptions, results from SOS programming can be leveraged to reduce the synthesis problem  in Theorem~\ref{thm:eventuality-safty} to a convex semidefinite optimization.

\subsection{System Dynamics and Semi-Algebraic Sets}
\label{sec:semialgebraic_formulation}
        The dynamics of both the evader and the pursuer are assumed to be linear in the control input, satisfying $f_p = f_e = \mathbf{0}_n$ and $g_e = g_p = \mathbf{I}_n$ with a norm constraint imposed on the control input $u$. Thus, the full system dynamics are defined as:
\begin{equation}
\label{eq:combined_system_simple}
    \bm{\dot{\mathbf{x}}} = \begin{bmatrix}
        \mathbf{u} \\ \mathbf{w}
    \end{bmatrix}, \qquad \| \mathbf{u} \|_\infty \leq u_{\text{max}}.
\end{equation}
This formulation can be approximately transformed into the Dubins car model via a suitable change of parameters, as demonstrated in~\cite{moyalan2024synthesizing}. This renders the problem to a game where the control actions $\mathbf{u}$ and $\mathbf{w}$ directly affect the positions of the evader and the pursuer. To define the initial, unsafe, and target sets that will be considered in this problem, we use the following functions:
\begin{equation}
    \begin{aligned}
    h_{\mathcal{X}e}(\mathbf{x}) &= x_1^2 + x_2^2 - R^2, \\
    h_{\mathcal{X}p}(\mathbf{x}) &= x_3^2 + x_4^2 - R^2, \\
    h_{ie}(\mathbf{x}) &= (x_1 - x_{ie_1})^2 + (x_2 - x_{ie_2})^2 - R_{ie}^2, \\
    h_{ip}(\mathbf{x}) &= (x_3 - x_{ip_1})^2 + (x_4 - x_{ip_2})^2 - R_{ip}^2, \\
    h_a(\mathbf{x}) &= (x_1 - x_3)^2 + (x_2-x_4)^2 - R_a^2, \\
    h_{re}(\mathbf{x}) &= (x_1 - x_{r_1})^2 + (x_2 - x_{r_2})^2 - R_r^2
    \label{eq:h}
    \end{aligned}
\end{equation}
The functions in \eqref{eq:h} define circular regions centered at different points in the state space. Using these functions, we define the sets in consideration as:
\begin{equation}
\begin{aligned}
    \mathcal{X} &= \{ \mathbf{x}: (h_{\mathcal{X}e} \leq 0 \vee h_{re} \leq 0) \wedge (h_{\mathcal{X}p} \leq 0 \} \\
    \mathcal{X}_i &= \{ \mathbf{x}: (h_{ie} < 0) \wedge (h_{ip} < 0) \} \\
    \mathcal{X}_a &= \{ \mathbf{x}: (h_{\mathcal{X}e} \leq 0) \wedge (h_{\mathcal{X}p} \leq 0) \wedge (h_a \leq 0) \} \\
    \mathcal{X}_r &= \{ \mathbf{x}: (h_{\mathcal{X}e} \geq 0) \wedge (h_{re} \leq 0) \wedge (h_{\mathcal{X}p} \leq 0) \} \\
\end{aligned}
\label{eq:sets_X}
\end{equation}
The sets \eqref{eq:sets_X} satisfy the conditions of Theorem \ref{thm:eventuality-safty}: $\mathcal{X}_i$, $\mathcal{X}_a$, and $\mathcal{X}_r$ are contained in $\mathcal{X}$; $\mathcal{X}$ is bounded; $\mathcal{X}_i$ is open; and $\mathcal{X}_i$, $\mathcal{X}_a$, and $\mathcal{X}_r$ are mutually disjoint. Therefore, we have
\begin{equation}
\begin{aligned}
    &\mathcal{X}_c = \operatorname{cl} ( \mathcal{X} \setminus (\mathcal{X}_r \cup \mathcal{X}_a) ) = \\
    &\quad \{ (h_{\mathcal{X}e} \leq 0) \wedge (h_{\mathcal{X}p} \leq 0) \wedge (h_a \geq 0) \} \\
    &\operatorname{cl}(\partial \mathcal{X} \setminus \partial \mathcal{X}_r) \cup \mathcal{X}_a = \\
    &\quad \{ (h_{\mathcal{X}e} = 0 \wedge h_{re} > 0) \vee (h_{\mathcal{X}p} = 0) \\
    &\quad \vee (h_a \leq 0 \wedge h_{\mathcal{X}e} \leq 0 \wedge h_{\mathcal{X}p} \leq 0) \}, \\
    &\operatorname{cl}(\mathcal{X} \setminus \mathcal{X}_r) = \{ (h_{\mathcal{X}e} \leq 0) \wedge (h_{\mathcal{X}p} \leq 0) \}.
\end{aligned}
\end{equation}

\subsection{Designing the Controller}
\label{sec:design_controller}
With the assumptions given in Section~\ref{sec:semialgebraic_formulation}, we obtain
\begin{equation*}
\begin{aligned}
    & \rhon (f +  g_w \mathbf{w}) + g_u \psin = 
    \begin{bmatrix}
        \rhon f_e + g_e \psin \\ \rhon f_p + \rhon g_p \mathbf{w}
    \end{bmatrix} = \begin{bmatrix} \psin \\ \rhon \mathbf{w} \end{bmatrix}.
\end{aligned} 
\end{equation*}
Thus, the divergence condition \eqref{eq:divergence_condition} can be restated as
\begin{small}
\begin{equation}
\begin{aligned}
    \nabla \cdot \begin{bmatrix} \psin \\ \rhon \mathbf{w} \end{bmatrix} > \alpha \nabla \log(V) \cdot \begin{bmatrix} \psin \\ \rhon \mathbf{w} \end{bmatrix} \quad \forall \mathbf{x} \in \operatorname{cl}(\mathcal{X} \setminus \mathcal{X}_r), \; \forall \mathbf{w} \in \mathcal{W}.
\end{aligned}
\label{eq:convex_rho_psi}
\end{equation}
\end{small}
Expanding the divergence and gradient operators, the inequality in \eqref{eq:convex_rho_psi} can be written as
\begin{equation}
    \nabla_e \cdot (\psin) + ( \nabla_p \rhon - \alpha \rhon \nabla_p \log(V) ) \cdot \mathbf{w} > \alpha \nabla_e \log(V) \cdot \psin.
\label{eq:convex_div_exp_2}
\end{equation}
where $\nabla_e$ and $\nabla_p$ denotes the divergence or gradient w.r.t. $\mathbf{x}_e$ and $\mathbf{x}_p$ respectively. 
 Define the set $\mathcal{P}$ as all $\mathbf{w} \in \mathbb{R}^{n_p}$ satisfying \eqref{eq:convex_div_exp_2}. Additionally, model the pursuer's movements with bounded velocity:
\begin{equation}
    \mathcal{W} = \left\{\mathbf{w} : 
    \begin{bmatrix}
        \mathbf{I} \\
        -\mathbf{I}
    \end{bmatrix}\mathbf{w}
    \leq
    w_{\max}\mathbf{1} \right\}.
\label{eq:W} 
\end{equation}
Robust weak eventuality and evasion are achieved if \eqref{eq:convex_div_exp_2} holds for all $\mathbf{w} \in  \mathcal{W}$, or, equivalently, iff $\mathcal{W} \subseteq \mathcal{P}$. From the generalized Farkas' lemma \cite{henrion1999control}, this inclusion is equivalent to the existence of a vector function $\mathbf{y}(\mathbf{x}) \geq \mathbf{0}$ such that:
\begin{equation}
\label{eq:dual_w}
\begin{aligned}
    \mathbf{y}^T \mathbf{N} &= -\nabla_p \rhon  + \alpha \rhon \nabla_p \log(V) \\ \mathbf{y}^T \mathbf{e} &< \nabla_e \cdot \psin - \alpha \nabla_e \log(V) \cdot \psin    
\end{aligned}
\end{equation}
where
\begin{equation}
\label{eq:Ne}
    \mathbf{N} = \begin{bmatrix}
        \mathbf{I} \\
        -\mathbf{I}
    \end{bmatrix}, \;
    \mathbf{e} = w_{\max} \mathbf{1}.
\end{equation}
These results leads to the following theorem, providing a sufficient condition for robust weak eventuality and evasion. 
\begin{theorem}
\label{thm:control}
Let the assumptions in Section~\ref{sec:semialgebraic_formulation} hold. If there exist functions $\psin:\Rx{n} \rightarrow \Rx{n}$, $\rhon:\Rx{n} \rightarrow \mathbb{R}$, and $\mathbf{y}: \Rx{n} \rightarrow \mathbb{R}$ such that
\begin{small}
\begin{flalign}
    &\rhon \geq 0 &&\forall \mathbf{x} \in \mathcal{X}_i \label{eq:thm_3_Xi} \\
    &\rhon < 0 &&\forall \mathbf{x} \in \operatorname{cl}(\partial\mathcal{X} \setminus \partial \mathcal{X}_r) \cup \mathcal{X}_a \label{eq:thm_3_cl} \\
    &\mathbf{y} \geq \mathbf{0} &&\forall \mathbf{x} \in \mathbb{R}^n \label{eq:thm_3_y} \\
    &V \mathbf{y}^\top \mathbf{N} = - V \nabla_p \rhon + \alpha \rhon \nabla_p V &&\forall \mathbf{x} \in \operatorname{cl}(\mathcal{X} \setminus \mathcal{X}_r) \label{eq:thm_3_yN} \\
    &V \mathbf{y}^\top \mathbf{e} < V \nabla_e \cdot \psin - \alpha \nabla_e V \cdot \psin &&\forall \mathbf{x} \in \operatorname{cl}(\mathcal{X} \setminus \mathcal{X}_r) \label{eq:thm_3_Vy} \\
    &| \psin | \leq u_{\max} \rhon &&\forall \mathbf{x} \in \mathcal{X}_c \label{eq:thm_3_umax}
\end{flalign}
\end{small}
Then, the control input $\mathbf{u} = \frac{\psin}{\rhon}$ enforces the \emph{weak eventuality and evasion} properties on \eqref{eq:combined_system_simple} in a \emph{robust manner} i.e. for almost all $\mathbf{x}_0 \in \mathcal{X}_i$, $\exists T \geq 0 \; \phi_T(\mathbf{x}_0) \in \mathcal{X}_r$ and $ \forall t \in [0,T] \; \phi_t(\mathbf{x}_0) \notin \mathcal{X}_a$ for all pursuers actions $\mathbf{w} \in \mathcal{W}$.
\end{theorem}

\subsection{SOS Implementation}
\label{sec:convex_relaxation_sos}
The conditions in Theorem~\ref{thm:control} can be formulated using Putinar’s Positivstellensatz theorem \cite{lasserre2009moments} as follows:
\begin{align*}
    \eqref{eq:thm_3_Xi}: \quad & \rhon + \sigma_{ie} h_{ie} + \sigma_{ip} h_{ip} &\in \Sigma[\mathbf{x}] \\
    \eqref{eq:thm_3_cl}: \quad & -\rhon + \lambda_{e} h_{\mathcal{X}e} -\sigma_{re} h_{re} &\in \Sigma[\mathbf{x}] \\
    & -\rhon + \lambda_{p} h_{\mathcal{X}p} &\in \Sigma[\mathbf{x}] \\
    & -\rhon + \sigma_a h_a + \sigma_{ae} h_{\mathcal{X}e} + \sigma_{ap} h_{\mathcal{X}p} &\in \Sigma[\mathbf{x}] \\
    \eqref{eq:thm_3_y}: \quad & \mathbf{y} &\in \Sigma[\mathbf{x}] \\
    \eqref{eq:thm_3_yN}: \quad & V \mathbf{y}^\top \mathbf{N} + V \nabla_p \rhon - \alpha \rhon \nabla_p V & = 0 \\
    \eqref{eq:thm_3_Vy}: \quad & -V \mathbf{y}^\top \mathbf{e} + V \nabla \cdot \psin - \alpha \nabla_e V \cdot \psin \\
    & + \sigma_{de} h_{\mathcal{X}e} + \sigma_{dp} h_{\mathcal{X}p} &\in \Sigma[\mathbf{x}] \\
    \eqref{eq:thm_3_umax}: \quad & \rhon u_{\max} \mathbf{1}_n - \psin + \sigma_{\mathbf{u}e1} h_{\mathcal{X}e} \\
    & + \sigma_{\mathbf{u}p1} h_{\mathcal{X}p} - \sigma_{\mathbf{u}a1} h_a &\in \Sigma[\mathbf{x}] \\
    & \rhon u_{\max} \mathbf{1}_n + \psin + \sigma_{\mathbf{u}e2} h_{\mathcal{X}e} \\
    & + \sigma_{\mathbf{u}p2} h_{\mathcal{X}p} - \sigma_{\mathbf{u}a2} h_a &\in \Sigma[\mathbf{x}] \\
    (\sigma): \quad & \sigma_{ie}, \sigma_{ip}, \sigma_{re}, \sigma_a, \sigma_{de}, \sigma_{dp} &\in \Sigma[\mathbf{x}] \\
    & \sigma_{\mathbf{u}e1}, \sigma_{\mathbf{u}e2}, \sigma_{\mathbf{u}p1}, \sigma_{\mathbf{u}p2}, \sigma_{\mathbf{u}a1}, \sigma_{\mathbf{u}a2} &\in \Sigma[\mathbf{x}]\\
    (\lambda): \quad &\lambda_e,\lambda_p & \in \mathbb{R}[\mathbf{x}]
\end{align*}
The $(\sigma)$ condition ensures the polynomial multipliers $\sigma$ are SOS, while $\lambda_e$ and $\lambda_e$ do not need to be SOS as they concern boundary conditions. Solving this feasibility problem gives us $\psin$ and $\rhon$, and the corresponding  control action $\mathbf{u} = \frac{\psin}{\rhon}$.

\begin{figure}[b]
    \vspace{-10pt}
    \centering
    \begin{minipage}{0.23\textwidth}
        \centering
        \includegraphics[width=1\textwidth]{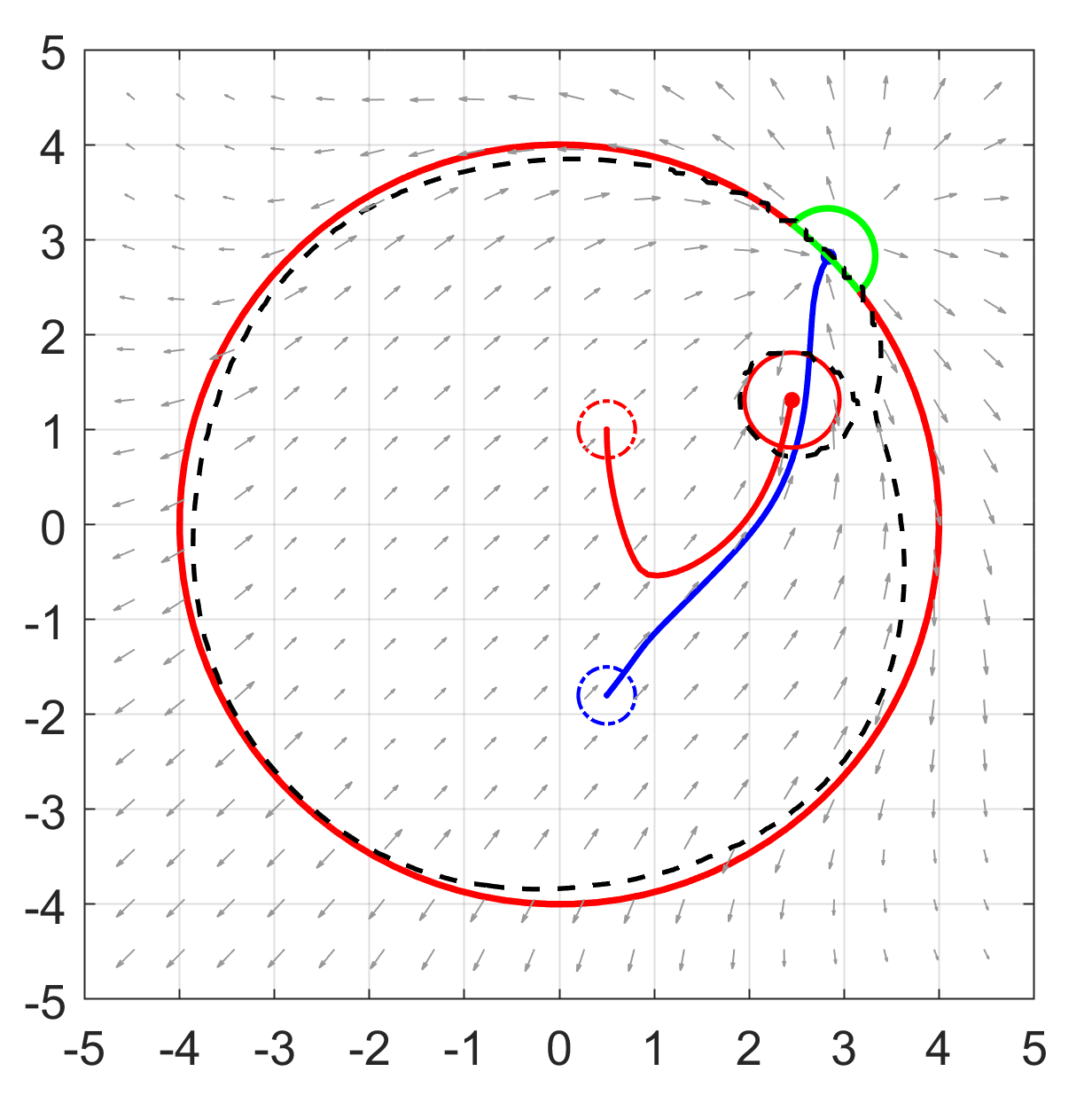}
    \end{minipage} 
    \begin{minipage}{0.23\textwidth}
        \centering
        \includegraphics[width=1\textwidth]{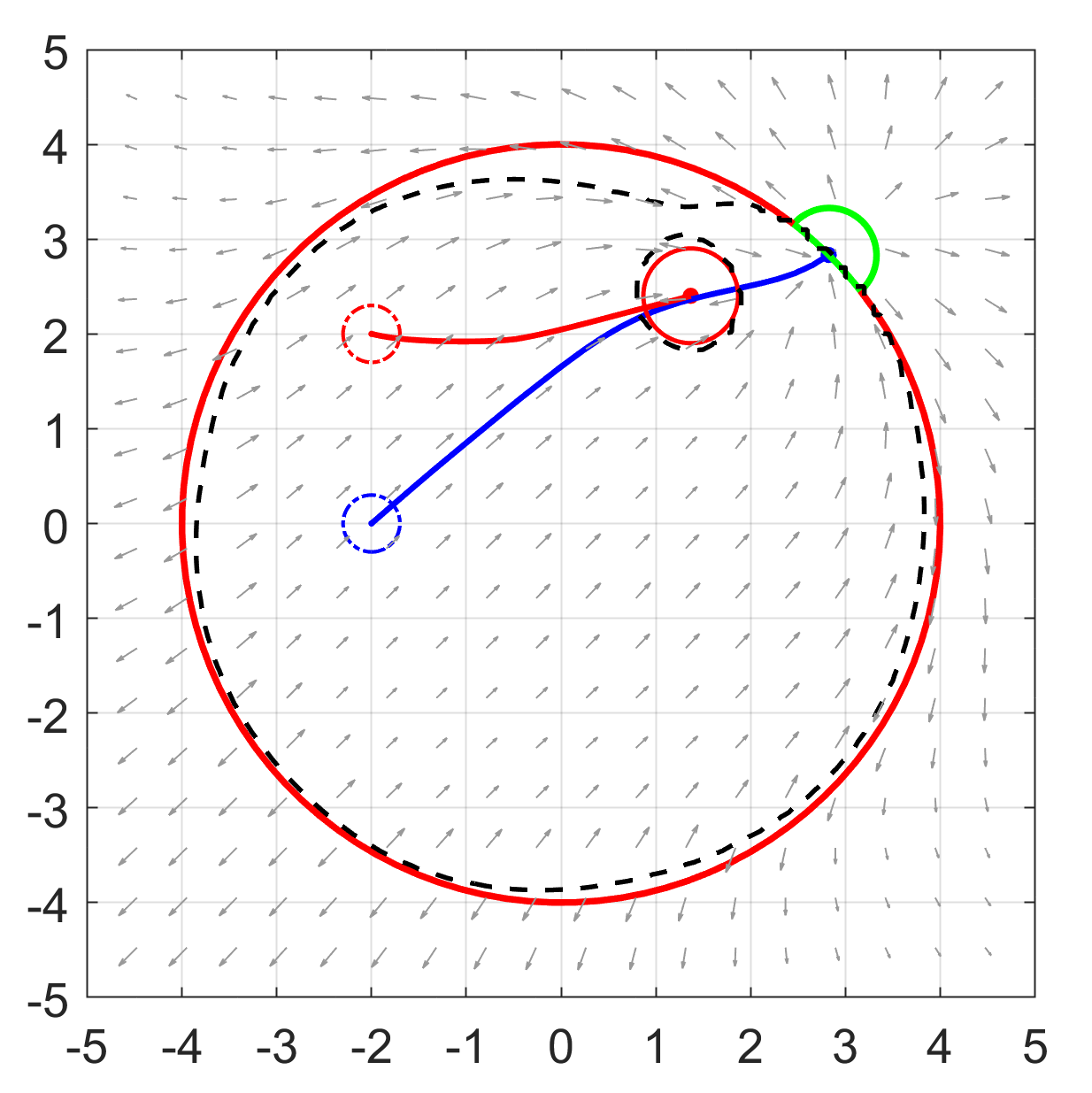}
    \end{minipage}
    \vspace{0pt}
    \caption{Final trajectories with tail-chasing (left) and go-to-middle (right) pursuers. The \textcolor{blue}{evader} reaching the \textcolor{green}{target set} while avoiding the \textcolor{red}{catch radius of the pursuer} and remaining in the \textcolor{red}{bounded region}, using the density function $\rho$. The dashed black line is the $\rho=0$ level set.}
    \label{fig:pe_simulation}
\end{figure}

\section{NUMERICAL EXAMPLE}  
\label{sec:numerical}

We evaluate the proposed method on a reach-avoid pursuit-evasion game with two pursuer strategies. Our code is available at \href{https://github.com/mbozdag08/pursuit-evasion-density}{pursuit-evasion-density}. The environment and region parameters in \eqref{eq:h} are defined as follows:
\begin{equation*}
    \begin{aligned}
    R &= 4, \quad R_{ie} = R_{ip} = 0.3, \quad R_a = R_r = 0.5, \\
    \mathbf{x}_r &= [R\cos(\pi/4), R\sin(\pi/4)].
    \end{aligned}
\end{equation*}
The target region lies on the boundary of $\mathcal{X}$ with radius $R$. The maximum control inputs are set as $u_{\max} = 0.015$ for the evader and $w_{\max} = 0.01$ for the pursuer, ensuring the evader has a feasible escape strategy. Polynomial degrees are chosen as $d_{\rhon} = d_{\psin} = 10$, $d_{V} = 2\alpha = 36$, with multiplier degrees $\deg(\sigma) = \deg(\lambda) = 6$. Two pursuer strategies are considered. In the first, the pursuer follows a \emph{tail-chasing} policy, always moving toward the evader at full speed. In the second, a \emph{go-to-middle} strategy is used, where the pursuer moves toward the midpoint between the evader and the center of the target arc. Initial positions are chosen as $\mathbf{x}_{ie}=[0.5, -1.8]$ and $\mathbf{x}_{ip}=[0.5, 1]$ for the tail-chasing case, and $\mathbf{x}_{ie}=[-2, 0]$ and $\mathbf{x}_{ip}=[-2, 2]$ for the go-to-middle case. The SOS program is implemented in MATLAB~\cite{MATLAB} using YALMIP~\cite{lofberg2004yalmip} and solved with MOSEK~\cite{mosek}.

\begin{figure}[t]
    \vspace{1pt}
    \centering
    \begin{minipage}{0.23\textwidth}
        \centering
        \includegraphics[width=1\textwidth]{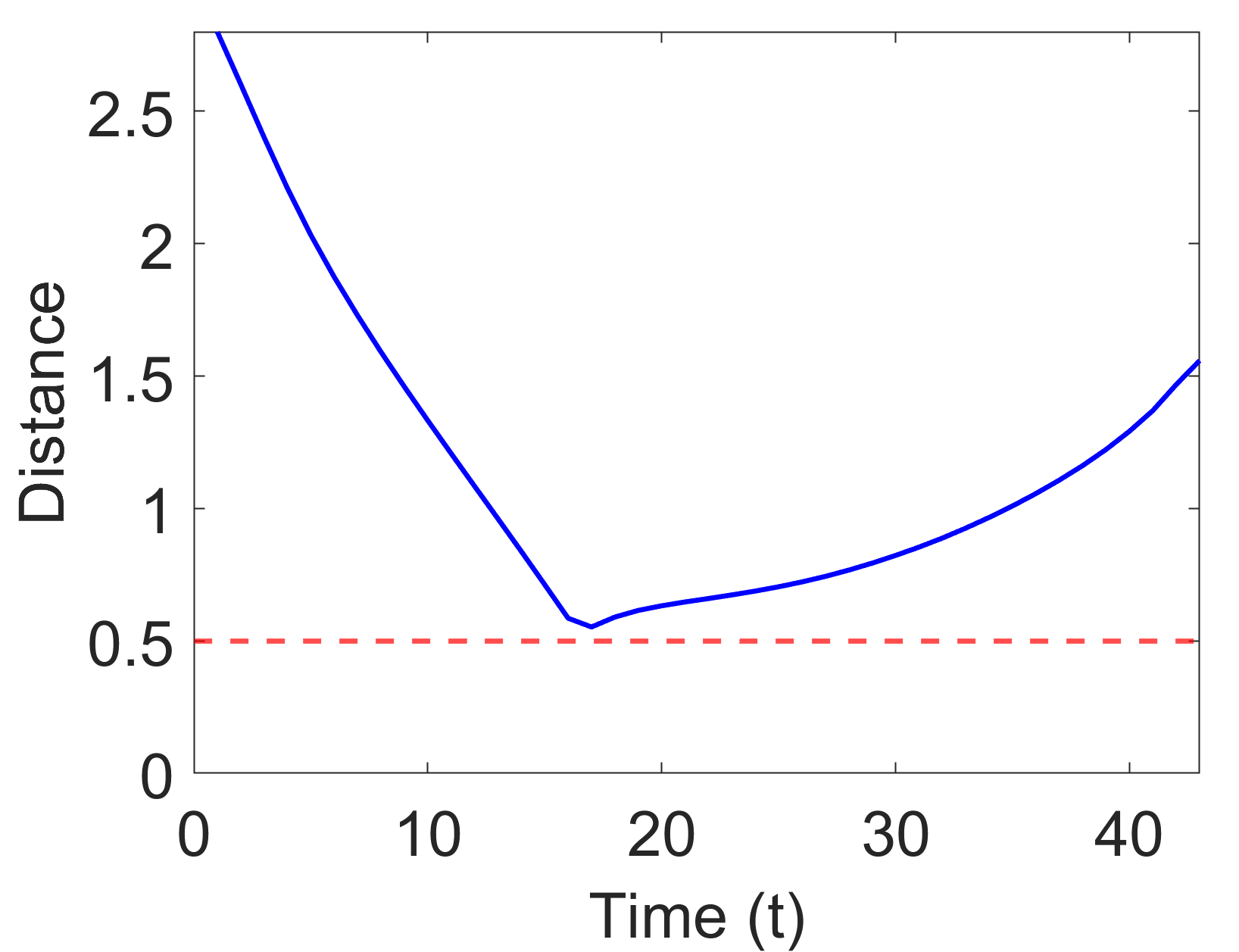}
    \end{minipage} 
    \begin{minipage}{0.23\textwidth}
        \centering
        \includegraphics[width=1\textwidth]{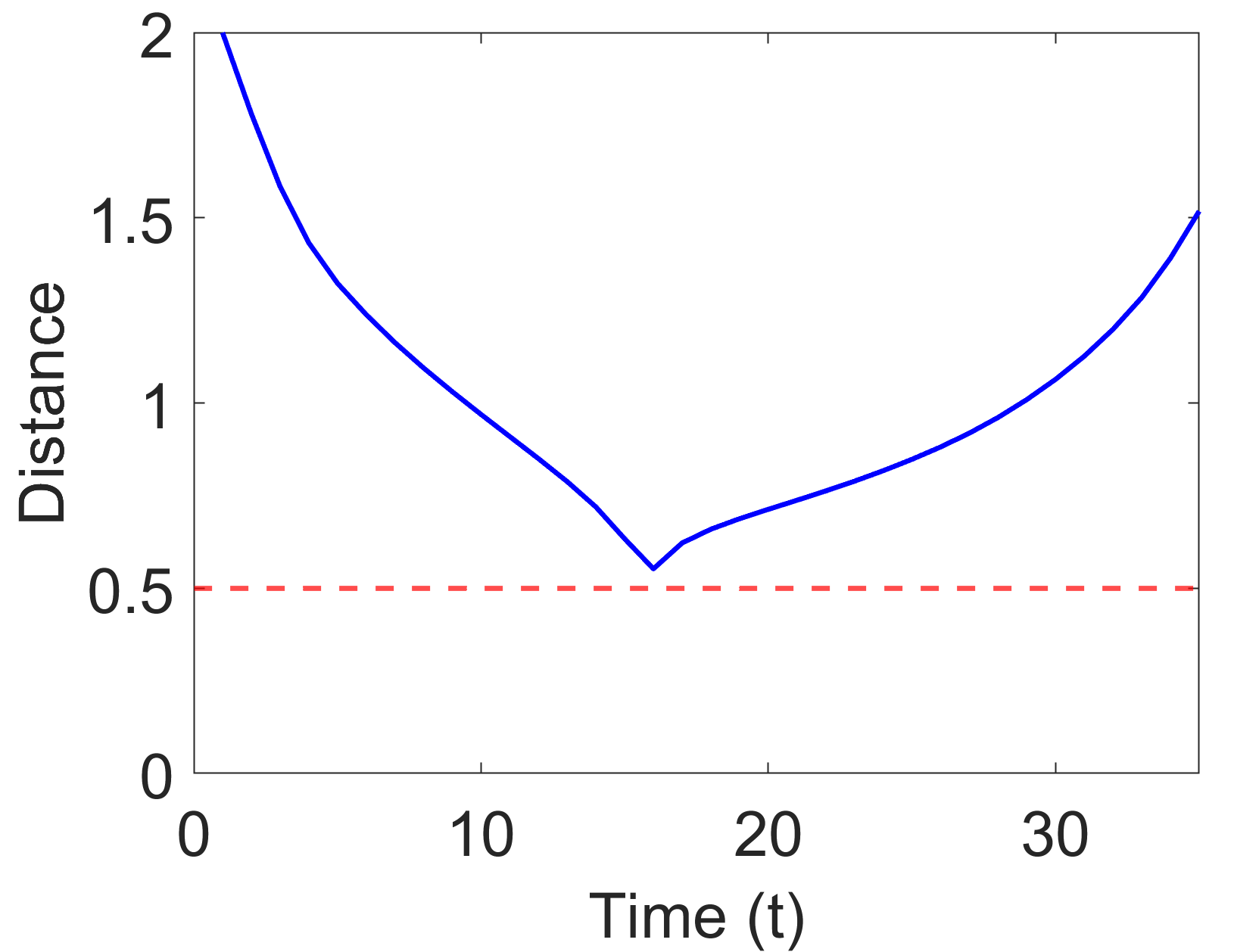}
    \end{minipage}
    \vspace{0pt}
    \caption{The distance between the evader and the pursuer over time for tail-chasing (left) and go-to-middle (right) pursuers. The dashed red line shows the catch-radius of the pursuer.}
    \label{fig:pe_distance}
    \vspace{-18pt}
\end{figure}

Figure~\ref{fig:pe_simulation} shows final-time snapshots of both simulations. The environment $\mathcal{X}$ and the pursuer’s unsafe zone $\mathcal{X}_a$ are shown as solid red circles. Dotted circles mark initial sets, and trajectories are shown as solid lines. In both scenarios, the evader (blue trajectory) reaches the crescent-shaped target region $\mathcal{X}_r$ (green) while staying outside the unsafe region (red). Figure~\ref{fig:pe_distance} shows the Euclidean distance between the evader and the pursuer over time. In both cases, the distance remains above the capture radius, validating safety throughout the simulations.

\vspace{-10pt}

\section{CONCLUSION}
This letter uses density functions to solve a reach-avoid pursuit-evasion game from a safe control perspective. Our approach offers a computationally efficient alternative to the classical HJI framework of differential game theory. The implemented SOS formulation yields a valid control action that satisfies the weak eventuality and evasion constraints of the game. Future work involves extending the density function approach to other variations of the game and with more complicated, unknown pursuer dynamics.
\label{sec:conclusion}


\bibliographystyle{IEEEtran}
\bibliography{references}

\end{document}

%% file: formulas.tex
\newcommand{\zero}
{
\mathbf{0}
}

\newcommand*{\Rx}[1]
{
\mathbb{R}^{#1}
}

\newcommand{\rhon}
{
\bar{\rho}
}

\newcommand{\psin}
{
\bm{\bar{\uppsi}}
}